\documentclass[eujc]{academic}

\def\Value#1{v(#1)}

\begin{document}
\shortauthor{P.~Vojt\v echovsk\'y}
\shorttitle{Reconstruction of Multiplication Tables}
\title{Reconstruction of Group Multiplication Tables by Quadrangle
Criterion}
\author{Petr Vojt\v echovsk\'y}
\address{Department of Mathematics, Iowa State University, 400 Carver Hall,
Ames, IA, 50011, U.S.A.}

\maketitle

\abstract{For $n>3$, every $n\times n$ partial Cayley matrix with at
most $n-1$ holes can be reconstructed by quadrangle criterion. Moreover, the
holes can be filled in given order. Without additional assumptions, this is
the best possible result. Reconstruction of other types of multiplication
tables is discussed.}

\section{Introduction}

\noindent Let us get started by explaining exactly what we mean by
reconstruction of multiplication tables. There are at least two approaches in
the literature, so it is not out of place to introduce the basic definitions
here.

Let $M=(m_{i,j})$ be an $n\times n$ matrix. It is important to distinguish the
entries of $M$ from the cells they occupy. For that matter, if $c=(i,\,j)$ is
the cell formed by the intersection of row $i$ and column $j$, let $\Value{c}$
be the entry in $c$, namely $m_{i,j}$.

The quadruple $(a,\,b,\,c,\,d)$ of cells is called \emph{quadrangle} if $a$,
$b$, $c$, $d$ are (all four) corners of a non-degenerate rectangular block
such that $a$, $c$ lie on one of the diagonals of the block. A block is
\emph{non-degenerate} if it has at least two rows and two columns. There are
$8$ ways to write down every quadrangle, and we will identify them.

A matrix $M$ is said to satisfy the \emph{quadrangle criterion} if
$\Value{c_4}=\Value{d_4}$ whenever $(c_1,\,c_2,\,c_3,\,c_4)$ and
$(d_1,\,d_2,\,d_3,\,d_4)$ are two quadrangles satisfying
$\Value{c_i}=\Value{d_i}$ for $i=1$, $2$, $3$. This criterion was introduced
by M.~Frolov \citep{Frolov}, as remarked by D\'enes and Keedwell \citep[p.\
19]{DenesKeedwell}.

Following \citep{DOP}, we say that $M$ is a \emph{Cayley matrix} if it is a
Latin square satisfying the quadrangle criterion.

An $n\times n$ matrix $M$ with a headline and sideline is a \emph{$($group$)$
multiplication table}, or \emph{Cayley table}, if there is a group
$(G,\,\cdot)$ and two enumerations $g_1$, $\dots$, $g_n$ and $h_1$, $\dots$,
$h_n$ of its elements such that the rows of $M$ are labeled by $g_1$, $\dots$,
$g_n$, the columns of $M$ by $h_1$, $\dots$, $h_n$ (in this order) and
\begin{equation}\label{Eq:CT}
    m_{i,\,j}=g_i\cdot h_j
\end{equation}
holds for every $i$, $j$, $1\le i$, $j\le n$.

The relation between Cayley tables and Cayley matrices is well-known (cf.
\citep[Theorem 1.2.1]{DenesKeedwell}). Every Cayley table gives rise to a
Cayley matrix when the sideline and headline are deleted. Conversely, given a
Cayley matrix $M$, any row and any column can be chosen as a sideline and
headline, respectively, to turn $M$ into a Cayley table. In other words, a
matrix $M$ is a Cayley matrix if and only if (\ref{Eq:CT}) holds for some
group $(G,\,\cdot)$ and two enumerations $g_1$, $\dots$, $g_n$ and $h_1$,
$\dots$, $h_n$ of elements of $G$.

A Cayley matrix $M$ will be called \emph{balanced} here, if $g_i=h_i$ for
every $i$, $1\le i\le n$. (We prefer the adjective balanced to
\emph{symmetrical} because a balanced Cayley matrix $M$ is symmetrical as a
matrix if and only if it is associated with a commutative group.) Similarly,
we can speak of \emph{balanced Cayley tables}.

There are other \emph{types} of multiplication tables. For instance,
Zassenhaus introduced normal Cayley tables, and Tamari defined generalized
normal Cayley tables (cf. \citep[p.\ 21]{DenesKeedwell}).

We are concerned with reconstructions of partial Cayley tables and matrices.
The process of reconstruction is not well-defined unless we specify:
\begin{enumerate}
\item[(a)] \ the type of multiplication table,
\item[(b)] \ the data available for the reconstruction,
\item[(c)] \ the method of reconstruction.
\end{enumerate}
The procedure then goes as follows. Let $G$ be a group, and $M$ one of its
multiplication tables of type $t$. Let us delete a few entries in $M$. The
resulting partial table $P$ will be referred to as a \emph{partial
multiplication table}, and the empty cells will be called \emph{holes}. Our
goal is to fill in the holes of $P$ using only the allowed data and methods so
that $P$ turns into a multiplication table of type $t$ again. When this
process, called \emph{reconstruction}, always yields $M$, we say that $M$ is
\emph{reconstructable} from $P$.

Here are a few comments on items $(a)$, $(b)$ and $(c)$. Naturally, we assume
that $P$ is part of the data available for reconstruction. However, when $M$
has no headline and sideline, it can happen that $P$ does not contain all
elements of $G$. Then, strictly speaking, $M$ cannot be reconstructed from $P$
unless we include the elements of $G$ as part of the data. Some authors take
this for granted, of course. As far as the methods are concerned, one can
always use the most general of them---the \emph{brute force} method. (Fill in
$P$ at random. Check whether you have obtained a multiplication table of type
$t$. Do it in all possible ways.) Needless to say, such an approach is merely
of theoretical interest, and it is therefore essential to specify the methods.
Also, the parallelism of the reconstruction is of practical importance.

The crucial question is: \emph{Given a group multiplication table $M$ of type
$t$, how many holes can there be in $P$ so that $M$ is reconstructable from
$P$ using only the allowed data and methods}?

Note that reconstructable means \emph{uniquely} reconstructable, by our
definition.

D\'enes proved \citep[Section 3.2]{DenesKeedwell} that, with two exceptions for
$n=4$ and $6$, every Cayley matrix with at most $2n-1$ holes can be
reconstructed provided all the group elements are known. He only used the
quadrangle criterion and the fact that every Cayley matrix is a Latin square,
but he was not interested in the order in which the holes can be filled. His
proof was made more precise by Frische \citep{Frische}.

Dr\'apal investigated Cayley tables and proved \citep{Drapal} that every such
table is reconstructable if $n\ge 51$ and there are not more than about $6n$
holes (see \citep{Drapal} for precise statement). The case when $n$ is prime
was resolved in \citep{PetrOlomouc}.

In all three situations, the estimate on the number of holes cannot be improved
in general.

The author is not aware of any result concerning the reconstruction of
balanced Cayley matrices. Apparently, the problem of reconstruction for Cayley
tables is equivalent to that of balanced Cayley tables.

It is easy to see that the reconstructions of Cayley matrices, balanced Cayley
matrices, and Cayley tables pose three distinct problems. We can illustrate
this already for $n=3$. To avoid trivialities, assume that the group elements
are known and call them $a$, $b$, $c$. The partial Cayley matrix $1$ in Figure
\ref{Fg:ThreeProblems} cannot be reconstructed. There are two possibilities to
complete $1$ into a Cayley matrix ($2$ and $3$). However, when we know that $1$
is a partial balanced Cayley matrix, it must be symmetrical, and therefore $2$
is the only solution. Matrix $4$ cannot be reconstructed as a balanced Cayley
matrix (both $2$ and $5$ are solutions). However, when we label the rows and
columns of $4$ as in $6$, say, the element $a$ becomes the neutral element, and
thus $6$ can only be completed into $3$ as a Cayley table with the given
headline and sideline.

\begin{figure}
\begin{displaymath}
\begin{array}{lll}
    \begin{array}{c|ccc}
        1   &       &       &       \\
        \hline
            &   a   &   b   &   c   \\
            &       &       &       \\
            &       &       &       \\
    \end{array}
    \quad\quad
    &
    \begin{array}{c|ccc}
        2   &       &       &       \\
        \hline
            &   a   &   b   &   c   \\
            &   b   &   c   &   a   \\
            &   c   &   a   &   b   \\
    \end{array}
    \quad\quad
    &
    \begin{array}{c|ccc}
        3   &       &       &       \\
        \hline
            &   a   &   b   &   c   \\
            &   c   &   a   &   b   \\
            &   b   &   c   &   a   \\
    \end{array}

    \\
    \\
    \begin{array}{c|ccc}
        4   &       &       &       \\
        \hline
            &   a   &       &       \\
            &       &       &       \\
            &       &       &       \\
    \end{array}
    \quad\quad
    &
    \begin{array}{c|ccc}
        5   &       &       &       \\
        \hline
            &   a   &   c   &   b   \\
            &   c   &   b   &   a   \\
            &   b   &   a   &   c   \\
    \end{array}
    \quad\quad
    &
    \begin{array}{c|ccc}
        6   &   a   &   b   &   c   \\
        \hline
        a   &   a   &       &       \\
        c   &       &       &       \\
        b   &       &       &       \\
    \end{array}
\end{array}
\end{displaymath}
\caption{Partial Cayley Matrices and Tables}
\label{Fg:ThreeProblems}
\end{figure}

In this short note we prove:

\begin{theorem}\label{Th:Main}
Let $n>3$. Every Cayley matrix of order $n$ with at most $n-1$ holes can be
reconstructed by quadrangle criterion. Moreover, the order in which the holes
are to be filled can be chosen in advance.
\end{theorem}

We also argue that this is, in a sense, the best possible result.

\section{The Reconstruction}

\noindent Let $M$ be a Cayley matrix associated with some $n$-element group
$G$, and let $P$ be a part of $M$. The following obvious Lemma tells us how to
apply the quadrangle criterion during reconstruction.

\begin{lemma}\label{Lm:QC}
Let $(c_1,\,c_2,\,c_3,\,c_4)$, $(d_1,\,d_2,\,d_3,\,d_4)$ be two quadrangles in
$P$ and assume that $\Value{c_i}=\Value{d_i}$ holds for $i=1$, $2$, $3$, that
$c_4$ is not a hole, and that $d_4$ \emph{is} a hole. Than the hole $d_4$ must
be filled with $\Value{c_4}$ in order to complete $P$ into a Cayley matrix.
\end{lemma}

Therefore, we are done with Theorem \ref{Th:Main} as soon as we prove:

\begin{proposition}\label{Pr:Recon}
Let $n>3$, and assume that there are at most $n-1$ holes in $P$. Then for every
hole $d_4$ there are two quadrangles $(c_1,\,c_2,\,c_3,\,c_4)$ and
$(d_1,\,d_2,\,d_3,\,d_4)$ such that $\Value{c_i}=\Value{d_i}$ holds for every
$i=1$, $2$, $3$, and such that the only hole among the $c_i$'s and $d_i$'s is
$d_4$.
\end{proposition}

Equivalently, we can state Proposition \ref{Pr:Recon} as follows:

\begin{proposition}\label{Pr:Group}
Assume that $G$ is a group of order $n>3$, and let $T\subseteq G\times G$ be
of cardinality at most $n-1$. Then for every tuple $(g_1,\,h_1)\in T$ there are
elements $g_2$, $h_2$, $g_1'$, $g_2'$, $h_1'$, $h_2'$ of $G$ such that
\begin{enumerate}
\item[(i)] \ $g_1\ne g_2$, $h_1\ne h_2$, $g_1'\ne g_2'$, $h_1'\ne h_2'$,
\item[(ii)] \ $g_i\cdot h_j=g_i'\cdot h_j'$ for every $i$, $j\in\{1,\,2\}$,
\item[(iii)] \ $\{(g_i,\,h_j),\, (g_i',\,h_j');\; 1\le i,\,j\le 2\}
    \cap T=\{(g_1,\,h_1)\}$.
\end{enumerate}
\end{proposition}

We now prove Proposition \ref{Pr:Recon}. Suppose that $n>3$ and that there are
at most $n-1$ holes in $P$.

\begin{lemma}\label{Lm:D}
Assume that $(a$, $b$, $c$, $d)$ is a quadrangle satisfying
\begin{enumerate}
    \item[$(C1)$] \ $d$ is a hole,
    \item[$(C2)$] \ $a$, $b$, $c$ are not holes,
    \item[$(C3)$] \ $\Value{a}\ne\Value{c}$, $\Value{b}\ne\Value{d}$ in $M$.
\end{enumerate}
Then the hole $d$ can be filled by quadrangle criterion.
\end{lemma}
\begin{proof}
There are exactly $n$ quadrangles $Q_i=(a_i$, $b_i$, $c_i$, $d_i)$ in $M$ with
$\Value{a_i}=\Value{a}$, $\Value{b_i}=\Value{b}$, $\Value{c_i}=\Value{c}$, and
$\Value{d_i}=\Value{d}$, for $i=1$, $\dots$, $n$. Because $\Value{a}$,
$\Value{b}$, $\Value{c}$, $\Value{d}$ are four different elements of $G$, no
two quadrangles $Q_i$, $Q_j$ have a corner in common. Since there are at most
$n-1$ holes in $P$, one of the quadrangles $Q_i$ is complete, say $Q_k$. Apply
Lemma \ref{Lm:QC} to $Q_k$ and $(a$, $b$, $c$, $d)$.
\end{proof}

Pick a hole $d$ in $P$. Without loss of generality, we may assume that
$d=(n,\, n)$. We try to find a quadrangle $(a$, $b$, $c$, $d)$ satisfying the
assumptions of Lemma \ref{Lm:D}. As we will see later, this is possible
whenever $n>4$. The case $n=4$ requires special treatment.

From now on, let all quadrangles $Q=(a$, $b$, $c$, $d)$ be written in such a
way that $a$ is the bottom-left corner, $b$ the top-left corner, and $c$ the
top-right corner of $Q$. Define
\begin{eqnarray*}
    T&=&\{c;\;c \textrm{\ a hole in $P$}\},\\
    T_0&=&\{c\in T;\; c=(i,\,j),\, 1\le i,\,j<n\},\\
    T_x&=&\{c\in T;\; c=(i,\,n),\, 1\le i<n\},\\
    T_y&=&\{c\in T;\; c=(n,\,j),\, 1\le j<n\},
\end{eqnarray*}
and let $t=|T|$, $t_0=|T_0|$, $t_x=|T_x|$, $t_y=|T_y|$. We have
\begin{equation}\label{Eq:Crucial}
    t_0+t_x+t_y+1=t\le n-1,
\end{equation}
because the sets $T_0$, $T_x$ and $T_y$ are disjoint and $d$ does not belong
to $T_0\cup T_x\cup T_y$.

Given $a=(n,\,j)$ with $j<n$ there are either $n-2$ or $n-3$ quadrangles
$(a,\,?,\,?,\,d)$ satisfying $(C1)$ and $(C3)$. (There are $n-2$ of them if
and only if there is a quadrangle $(a,\,b,\,c,\,d)$ with $\Value{a}=\Value{c}$,
$\Value{b}=\Value{d}$.) Therefore, there are at least $(n-1)(n-3)$ quadrangles
$(?,\,?,\,?,\,d)$ satisfying $(C1)$ and $(C3)$. (This estimate cannot be
improved in general. To see this, consider the standard Cayley table of any
cyclic group $C_n$ of odd order $n>1$.)

\begin{lemma}\label{Lm:Aux}
If $t_x+t_y\le 1$, there is at least one quadrangle $(?,\,?,\,?,\,d)$
satisfying $(C1)$, $(C2)$ and $(C3)$.
\end{lemma}
\begin{proof}
There are at least $(n-1)(n-3)$ quadrangles $(?,\,?,\,?,\,d)$ satisfying
$(C1)$ and $(C3)$. Every hole from $T_x$ affects at most $n-2$ of them, and so
does every hole from $T_y$. Every hole from $T_0$ affects at most one such
quadrangle. Thus, there are at least
\begin{displaymath}
    \tau=(n-1)(n-3)-t_0-(t_x+t_y)(n-2)
\end{displaymath}
quadrangles satisfying $(C1)$, $(C2)$, $(C3)$.

When $t_x+t_y=0$, we have $\tau\ge(n-1)(n-3) - (n-2) = n^2-5n+5>0$, for $n>3$.

Similarly, when $t_x+t_y=1$, we have $t_0\le n-3$, and consequently $\tau\ge
(n-1)(n-3)-(n-3)-(n-2)=n^2-6n+8$. This is positive when $n>4$.

Without loss of generality, assume that $t_x=1$, $t_y=0$, $n=4$. Delete the
unique row $i<n$ of $M$ for which $(i,\,n)$ is a hole to obtain an $(n-1)\times
n$ block $B$. For every $a=(n,\,j)$ with $j<n$, there is at least one
quadrangle $(a,\,b,\,c,\,d)$ such that $b,\,c\in B$ and
$\Value{a}\ne\Value{c}$. Moreover, when $\Value{n,\,j}=\Value{i,\,n}$, there is
another such quadrangle. Hence, there are at least $3+1=4$ such quadrangles.
At most $2$ of them satisfy $\Value{b}=\Value{d}$ because $\Value{d}$ appears
$3$ times in $B$. Hence, there are at least $4-2=2$ quadrangles
$(?,\,?,\,?,\,d)$ satisfying $(C1)$ and $(C3)$. Since $t_0\le 1$, we are done.
\end{proof}

When $t_0\ge n-3$, we have $t_x+t_y\le 1$, and Lemma \ref{Lm:Aux} applies.

When $0<t_0<n-3$, we have $t_x+t_y\le n-3$. Then $\tau$ from the proof of
Lemma \ref{Lm:Aux} is greater than or equal to $n-3-t_0>0$.

Finally, assume that $t_0=0$. We could change our point of view and conclude
that at least one hole of $P$ can be filled but, remember, we want to fill $d$.

Let $t_y<n-3$. There is at least one full cell $c$ in the $n$th column of $P$.
Therefore, we have at least $3$ quadrangles $(?,\,?,\,c,\,d)$ satisfying
$(C1)$ and $(C2)$. The condition $(C3)$ can exclude at most $2$ of them.
Similarly when $t_x<n-3$.

When both $t_x$ and $t_y$ are bigger than or equal to $n-3$,
(\ref{Eq:Crucial}) implies that $2(n-3)\le t_x+t_y\le n-2$. This means that
$n=4$, $t_x$, $t_y\ge 1$, $t_0=0$. It is enough to solve the case $t_x=t_y=1$.
Let $c_1$, $c_2$ be the two full cells in the $n$th column of $P$. Similarly,
introduce $a_1$, $a_2$ in the $n$th row of $P$. Let $B$ be the $3\times 3$
top-left block of $M$. The value $\Value{d}$ appears in $B$. Pick $k$, $l$
such that $1\le k$, $l\le 3$ and $\Value{k,\,l}=\Value{d}$. Then there are
$i$, $j\in\{1,\,2\}$ such that $\Value{c_i}$ is within the $k$th row of $B$
and $\Value{b_j}$ within the $l$th column of $B$. Because there is no hole in
$B$, we have found two quadrangles satisfying the assumptions of Lemma
\ref{Lm:QC}.

This completes the proof of Proposition \ref{Pr:Recon} and Theorem
\ref{Th:Main}.

\section{Discussion and Acknowledgements}

\noindent The bound $n>3$ cannot be improved. Consider the Cayley matrix
\begin{displaymath}
    \begin{array}{ccc}
        0&1&2\\
        \framebox{$1$}&2&0\\
        2&\framebox{$0$}&1
    \end{array}
\end{displaymath}
of $C_3$, and observe that none of the holes (framed cells) can be filled by
quadrangle criterion.

Lemma \ref{Lm:D} cannot be applied to every Cayley matrix of order $4$. Look,
for example, at the Cayley matrix
\begin{displaymath}
    \begin{array}{cccc}
        0&1&2&3\\
        1&2&3&\framebox{$0$}\\
        2&3&0&1\\
        3&\framebox{$0$}&1&\framebox{$2$}
    \end{array}
\end{displaymath}
of $C_4$.

Because there are at most $n-1$ holes in $P$, we do not need to add the names
of elements of $G$ as part of the data available for reconstruction. However,
even with this data added, the bound $n-1$ cannot be improved to $n$. Consider
a Cayley matrix $M$ with one row deleted. Then it is impossible to reconstruct
$M$ just by quadrangle criterion. It would be interesting to know whether this
is the only pathological situation for $n>4$.

The author would like to thank one of the referees for his/her useful comments
on reconstruction in general, and for the formulation of Proposition
\ref{Pr:Recon} as a group-theoretical result (cf. Proposition \ref{Pr:Group}).

\bibliographystyle{plain}

\end{document}